\documentclass[11pt, a4paper]{amsart}
\usepackage[top=35truemm,bottom=35truemm,left=35truemm,right=35truemm]{geometry}
\usepackage{setspace} 
\usepackage{amsmath, amsfonts,amsthm,amssymb,mathrsfs,amscd}
\usepackage{enumerate}
\usepackage[dvipdfmx]{graphicx}
\usepackage{ascmac}
\usepackage[arrow,matrix]{xy}
\usepackage{color}
\usepackage{array}
\usepackage{pifont}
\usepackage{longtable}
\usepackage{here}
\usepackage{braket}
\usepackage{mathtools}

\newtheorem{theo}{Theorem}
\newtheorem*{theo*}{Theorem}

\newtheorem{cor}{Corollary}
\newtheorem{lem}{Lemma}

\providecommand{\norm}[1]{\left\lVert#1\right\rVert}
\providecommand{\abs}[1]{\left\lvert#1\right\rvert}
\newcommand{\bvec}[1]{\mbox{\boldmath $#1$}}
\DeclareMathOperator{\Sym}{Sym}
\DeclareMathOperator{\Sim}{Sim}
\DeclareMathOperator{\Isom}{Isom}
\DeclareMathOperator{\Hom}{Hom}

\def\disp{\displaystyle}

\def\B{\mathbb{B}}

\def\H{\mathbb{H}}
\def\R{\mathbb{R}}

\def\U{\mathbb{U}}

\def\({\left(}
\def\){\right)}
\def\<{\langle}
\def\>{\rangle}

\begin{document}

\title[ ]{Locally rigid right-angled Coxeter groups with Fuchsian ends in dimension 5}
\author{Tomoshige Yukita}
\address{Department of Mathematics, School of Education, Waseda University, Nishi-Waseda 1-6-1, Shinjuku, Tokyo 169-8050, Japan}
\email{yshigetomo@suou.waseda.jp}
%\keywords{Coxeter group; local rigidity; }
\date{}
\thanks{}

\begin{abstract} 
In this paper, we construct a right-angled $5$-polytope $P$ of finite volume 
such that all the right-angled Coxeter groups with Fuchsian ends obtained from $P$ are locally rigid. 
\end{abstract}

\maketitle
%%%%%%%%%%%%%%%%%%%%%%%%%%%%%%%%%%%%%%%%%%%%%%%%%%%%%%%%%%%%%%%%%%%%%%%%%%%%%%%%%%%%%%%%%%%%%%%%%

\section{Introduction}
Let $\mathbb{H}^{d}$ denote the hyperbolic $d$-space. 
A non-empty subset $P\subset{\overline{\mathbb{H}}^d}$ is called a \textit{$d$-dimensional hyperbolic polytope} (\textit{$d$-polytope} for short) 
if $P$ can be written as the intersection of finitely many closed half-spaces $H_1^-, \cdots, H_N^-$. 
A hyperbolic polytope $P=\bigcap{H_i^-}$ is called a {\textit{hyperbolic right-angled polytope}} if all of its dihedral angles are $0$ and $\frac{\pi}{2}$. 
If $P\subset{\overline{\mathbb{H}}^{d}}$ is a hyperbolic right-angled polytope, 
the set $S=\{r_1, \cdots, r_N\}$ of all the reflections with respect to the facets of $P$ 
generates a discrete group $\Gamma_P$ of $\Isom{(\H^d)}$. 
We call $\Gamma_P$ a \textit{$d$-dimensional hyperbolic right-angled Coxeter group} (\textit{$RACG$} for short). 

\vspace{2mm}
Suppose that $\Gamma$ is a discrete subgroup of $\Isom{(\H^d)}$ and $\rho_0$ is the inclusion map. 
Let us denote by $\Hom{(\Gamma, \Isom{(\H^d)})}$ the space of all the group homomorphisms from $\Gamma$ to $\Isom{(\H^d)}$ 
with topology of pointwise-convergence, called the \textit{representation space of $\Gamma$}. 
$\Gamma$ is said to be \textit{locally rigid} if every representation $\rho$ nearby $\rho_0$ is a conjugate of $\rho_0$ by $\Isom{(\H^d)}$. 
For $d\geq{4}$, Calabi \cite{Ca},Weil \cite{We}, Garland, and Raghunathan \cite{GR} showed that every lattice $\Gamma$ is locally rigid. 
Thus we are interested in the rigidity of discrete groups of infinite covolume. 
In order to study the rigidity for such a case, Kerckhoff and Storm introduced the notion of discrete groups \textit{with Fuchsian ends} \cite{KS1}. 
In particular, Kerckhoff, Storm, and Aougab studied the local rigidity of RACGs with Fuchsian ends obtained from the 24-cell group and 120-cell group (\cite{AS, KS1}). 

\vspace{2mm}
The goal of this work is the study of the local rigidity of RACGs with Fuchsian ends in dimension 5. 
In the context of \textit{convex cocompact} RACGs, the results by Kerckhoff and Storm \cite{KS2}, together with the colouring technique by Kolpakov and Slavich \cite{KoS}, showed that all the convex cocompact RACGs in dimension $d\geq{4}$ are locally rigid. 

\vspace{2mm}
In this paper, we construct a right-angled $5$-polytope $\mathcal{P}$ which is highly symmetric and 
study the local rigidity of the RACGs obtained from $\Gamma_{\mathcal{P}}$. 
More precisely, the $5$-polytope $\mathcal{P}$ is constructed by adding spheres to the Euclidean regular $4$-cube $\mathcal{C}$, so that the symmetry group $\Sym{(\mathcal{C})}$ can be identified index 2 subgroup of the symmetry group $\Sym{(\mathcal{P})}$. 
By using $\Sym{(\mathcal{P})}$, we can show that all the RACGs with Fuchsian ends obtained from $\Gamma_{\mathcal{P}}$ are locally rigid. 
In this way, we give the first explicit example of locally rigid discrete subgroups of $\Isom{(\H^5)}$ with Fuchsian ends. 

\vspace{2mm}
The paper is organized as follows. 
In Section 2, we give the necessary background about the right-angled Coxeter groups and its representation space. 
In Section 3, we construct a right-angled hyperbolic $5$-polytope $\mathcal{P}$ of finite volume and study its geometry and combinatorics. 
Then, in Section 4, we will show that all the RACGs with Fuchsian ends obtained from $\Gamma_{\mathcal{P}}$ are locally rigid. 

%%%%%%%%%%%%%%%%%%%%%%%%%%%%%%%%%%%%%%%%%%%%%%%%%%%%%%%%%%%%%%%%%%%%%%%%%%%%%%%%%%%%%%%%%%%%%%%%%

\section{Preriminalies}
The references here are \cite{Ra} and \cite{V}. 

%%%%%%%%%%%%%%%%%%%%%%%%%%%%%%%%%%%%%%%%%%%%%%%%%%%%%%%%%%%%%%%%%%%%%%%%%%%%%%%%%%%%%%%%%%%%%%%%%
\subsection{Hyperbolic $d$-space}
\textit{Lorentzian $(d+1)$-space $\R^{d,1}$} is a $(d+1)$-dimensional real vector space equipped with the \textit{Lorentzian inner product}
\begin{equation}
\<\bvec{x}, \bvec{y}\>=x_1y_1+\cdots{+}x_dy_d-x_{d+1}y_{d+1}. 
\end{equation}
Set $\norm{\bvec{x}}=\<\bvec{x}, \bvec{x}\>$. 
A vector $\bvec{x}$ of $\R^{d, 1}$ is said to be \textit{space-like, time-like}, and \textit{light-like} if 
$\norm{\bvec{x}}$ is positive, negative, and zero, respectively. 
Let us denote by $H_v$ the Lorentz orthogonal complement of $\bvec{v}$, that is, 
\begin{equation}
H_v=\Set{\bvec{x}\in{\R^{d,1}}|\<\bvec{x}, \bvec{v}\>=0}. 
\end{equation}
We write $H_v^-$ for the closed half-space with the outer normal vector $v$, that is, 
\begin{equation}
H_v^-=\Set{\bvec{x}\in{\R^{d,1}}|\<\bvec{x}, \bvec{v}\>\leq{0}}. 
\end{equation}

\vspace{2mm}
The hyperboloid 
$\H^d=\Set{\bvec{x}=(x_1, \cdots, x_{d+1})\in{\R^{d,1}} | \<\bvec{x}, \bvec{x}\>=-1, \ x_{d+1}>0}$ 
with the metric induced by $\R^{d,1}$ is a model of hyperbolic $d$-space, 
so called the \textit{hyperboloid model of hyperbolic $d$-space}. 
The isometry group $\Isom{(\H^d)}$ is identified with the positive Lorentz group $PO(d, 1)$, where 
\begin{equation}
PO(d, 1)=\Set{A\in{GL_{d+1}(\R)}|AJA^T=J}, \ J=\begin{pmatrix} I_d & \bvec{0} \\ \bvec{0} & -1 \end{pmatrix}. 
\end{equation}
For every space-like vector $\bvec{v}$, the intersection $\H^d\cap{H_v}$ is called a \textit{hyperplane} of $\H^d$. 
The intersection $\H^d\cap{H_v^-}$ is called the \textit{closed half-space with outer vector $\bvec{v}$}. 
By abuse of notation, we write $H_v$ and $H_v^-$ for the hyperplane $\H^d\cap{H_v}$ and the closed half-space $\H^d\cap{H_v^-}$, respectively. 
The boundary at infinity $\partial{\H^d}$ of $\H^d$ is identified with the set of all the light-like vectors of $x_{d+1}=1$. 

\vspace{2mm}
A non-empty subset $P\subset{\overline{\mathbb{H}}^d}$ is called a \textit{$d$-dimensional hyperbolic polytope} (\textit{$d$-polytope} for short) 
if $P$ can be written as the intersection of finitely many closed half-spaces. 
This means that $\disp P=\bigcap_{i=1}^N{H^-_i}$, where $H^-_i$ is the closed half-space of $\mathbb{H}^d$ bounded by the hyperplane 
$H_i$ with normal vector $\bvec{v}_i$ pointing outwards with respect to $P$. 
An intersection of $P$ and $H_i$ is called a \textit{facet of $P$}. 
Notice that every facet of $P$ is a $(d-1)$-polytope. 
The mutual positions of the bounding hyperplanes of $P$ are characterized by its Lorentzian inner products: 
\begin{enumerate}
\item[(i)] The hyperplanes $H_i$ and $H_j$ intersect $\Leftrightarrow$ $\abs{\<\bvec{v}_i, \bvec{v}_j\>}<1$. 
\item[(ii)] The hyperplanes $H_i$ and $H_j$ meet at a point at infinity $\Leftrightarrow$ $\abs{\<\bvec{v}_i, \bvec{v}_j\>}=1$. 
\item[(iii)] The hyperplanes $H_i$ and $H_j$ are disjoint $\Leftrightarrow$ $\abs{\<\bvec{v}_i, \bvec{v}_j\>}>1$.
\end{enumerate}
Suppose that $H_i\cap{H_j}\neq{\emptyset}$ in $\mathbb{H}^{d}$. 
Then, the {\textit {dihedral angle} $\theta_{ij}$} between $H_i$ and $H_j$ is defined by 
\begin{equation}
\cos{\theta_{ij}}=-\<\bvec{v}_i, \bvec{v}_j\>. 
\end{equation}
If $\overline{H_i}\cap{\overline{H_j}}\in{\overline{\mathbb{H}}^{d}}$ is a point on $\partial{\mathbb{H}^d}$, 
then the dihedral angle between $H_i$ and $H_j$ is defined to be zero. 

\vspace{2mm}
A hyperbolic polytope $\disp P=\bigcap_{i=1}^N{H_i^-}\subset{\overline{\mathbb{H}}^{d}}$ is called a {\textit{hyperbolic right-angled polytope}} if all of its dihedral angles are $0$ and $\frac{\pi}{2}$. 
If $P\subset{\overline{\mathbb{H}}^{d}}$ is a hyperbolic right-angled polytope, 
the set $S=\{r_1, \cdots, r_N\}$ of all the reflections with respect to the facets of $P$ 
generates a discrete group $\Gamma_P$ of $\Isom{(\H^d)}$. 
We call $\Gamma_P$ a \textit{$d$-dimensional hyperbolic right-angled Coxeter group} (\textit{$RACG$} for short). 
The pair $(\Gamma, S)$ is a \textit{Coxeter system}: 
\begin{equation}
\Gamma=\< \ S \ | \ r_1^2=\cdots =r_N^2=1, \ (r_ir_j)^2=1 \text{ if }H_i\text{ and }H_j\text{ intersect} \ \>. 
\end{equation}

\vspace{2mm}
A $d$-dimensional hyperbolic RACG $\Gamma$ has \textit{Fuchsian ends} if 
the convex core $C_\Gamma\subset{\H^d}$ is a $d$-manifold with nonempty totally geodesic boundaries 
and the quotient $C_\Gamma/\Gamma$ has finite volume. 
Let $F$ be a set of pairwise disjoint facets of a right-angled $5$-polytope $P$. 
We write $P_F$ for the right-angled $5$-polytope obtained from $P$ by removing the facets in $F$. 
The RACG associated with $P_F$ is denoted by $\Gamma_F$. 
Then, $\Gamma_F$ has Fuchsian ends. 
Moreover, the quotient $C_{\Gamma_F}/\Gamma_F$ is isometric to the original right-angled $d$-polytope $P$. 
The totally geodesic boundaries of $C_{\Gamma_F}/\Gamma_F$ are the elements of $F$. 
The RACG $\Gamma_F$ with Fuchsian ends is said to be \textit{obtained from $\Gamma_P$ by removing the facets in $F$}. 

\vspace{2mm}
The upper half-space 
$\U^d=\Set{(x_1, \cdots, x_d)\in{\R^d} | x_d>0}$ 
equipped with the metric $\frac{|dx|}{x_d}$ is a model of hyperbolic $d$-space, so-called the \textit{upper half-space model of hyperbolic $d$-space}. 
The boundary $\partial{\U^d}$ of $\U^d$ in the one-point compactification $\R^d\cup{\{\infty\}}$ of $\R^d$ is the \textit{boundary at infinity} of $\U^d$. 
Let us denote by $\zeta$ the stereographic projection of $\H^d$ onto the open unit ball $\B^d\subset{\R^d}$, that is, 
\begin{equation}
\zeta(\bvec{x})=\(\dfrac{x_1}{1+x_{d+1}}, \cdots, \dfrac{x_d}{1+x_{d+1}}\). 
\end{equation}
Define the standard transformation $\eta$ from $\B^d$ onto $\U^d$ by 
\begin{equation}
\eta(\bvec{y})=\dfrac{1}{\abs{\bvec{y}-\bvec{e}_d}^2}\(2y_1, \cdots, 2y_{d-1}, 1-\abs{\bvec{y}}^2\), 
\end{equation}
where $\bvec{e}_d=(0, \cdots, 0, 1)\in{\R^d}$. 
Then the composition $\Phi=\eta\circ{\zeta}$ is an isometry from $\H^d$ onto $\U^d$. 
By using $\Phi$, we can translate the terminology of $\H^d$ into that of $\U^d$. 
For an example, a subset $H\subset{\U^d}$ is said to be a hyperplane if 
it is either a Euclidean hemisphere or a half-plane in $\U^{d}$
orthogonal to $\R^{d-1}$. 
Thus every hyperplane $H$ of $\U^d$ is determined by the intersection $\overline{H}\cap{\R^{d-1}}$. 
We identify the hyperplanes and spheres of $\R^{d-1}$ with the hyperplanes of $\U^d$.

%%%%%%%%%%%%%%%%%%%%%%%%%%%%%%%%%%%%%%%%%%%%%%%%%%%%%%%%%%%%%%%%%%%%%%%%%%%%%%%%%%%%%%%%%%%%%%%%%
\subsection{Representation space and local rigidity}
Let $\Gamma=\Gamma_P$ be a $d$-dimensional RACG associated with $P$ 
and $\rho_0$ be the inclusion map from $\Gamma$ into $\Isom{(\H^d)}$. 
We denote by $\Hom{(\Gamma, \Isom{(\H^d)})}$ the space of all the group homomorphisms from $\Gamma$ to $\Isom{(\H^d)}$ 
with topology of pointwise-convergence, called the \textit{representation space of $\Gamma$}. 
The isometry group $\Isom{(\H^d)}$ acts on $\Hom{(\Gamma, \Isom{(\H^d)})}$ by conjugation. 
$\Gamma$ is said to be \textit{locally rigid} if 
there exists a neighborhood $U\subset{\Hom{(\Gamma, \Isom{(\H^d)})}}$ contained in the orbit of $\rho_0$, 
that is, for every $\rho\in{U}$, there exists an isometry $f\in{\Isom{(\H^d)}}$ such that $\rho=f\rho_0 f^{-1}$. 
The following lemma is a fundamental to study the rigidity of RACGs. 
\begin{lem}\label{lem:2-1}\cite[Corollary 2.4]{AS}
Let $S=\{r_1, \cdots, r_N\}$ be a generating set of $\Gamma$ satisfying that $\rho_0(r_i)$ is a reflection for $i=1, \cdots, N$. 
There exists an open neighborhood $U_0\subset{\Hom{(\Gamma, \Isom{(\H^d)})}}$ of $\rho_0$ such that 
if $\rho\in{U_0}$, then $\rho(r_i)$ is a reflection for $i=1, \cdots, N$. 
\end{lem}
We can assume, by shrinking $U_0$ if necessary, that $U_0$ is connected. 
For every $\rho\in{U_0}$, we write $H_i(\rho)$ for the hyperplane of $\H^d$ corresponding to the reflection $\rho(r_i)$. 
We choose an unit outer normal vector $\bvec{v}_i(\rho)$ of $H_i(\rho)$ so that the function $\bvec{v}_i:U_0\to{\R^{d, 1}}$ is continuous. 
Since $U_0$ is connected, such a choice is uniquely determined by $\rho_0$. 

\vspace{2mm}
We shall mention some known results on rigidity of RACGs. 
If $\Gamma$ has finite covolume and $d\geq{4}$, then $\Gamma$ is locally rigid by Calabi \cite{Ca}, Weil \cite{We}, Garland, and Raghunathan \cite{GR}. 
Therefore our aim is to study the local rigidity of RACGs of infinite covolume. 
Note that the Mostow rigidity theorem dose not imply the local rigidity of lattices, 
since the theorem only concerned with the \textit{discrete faithful representations} of lattices. 
Thus the local rigidity of lattices mean that every representation near by the inclusion map is always discrete faithful. 

Kerckhoff, Storm, and Aougab studied rigidity of RACGs with Fuchsian ends and showed that the followings: 
first, Aougab and Strom showed that the RACGs with Fuchsian ends obtained from the compact right-angled 120-cell group are locally rigid \cite{AS}. 
After that, Kerckhoff and Storm generalized the result to convex cocompact $d$-dimensional RACGs with Fuchsian ends for $d\geq{4}$ \cite{KS2}. 
Here, convex cocompact means that the orbifold with totally geodesic boundaries $C_\Gamma/\Gamma$ is compact. 
They also found an example of not locally rigid RACG which is not convex cocompact \cite{KS1}. 
From the not locally rigid example, Martelli and Riolo constructed the first example of hyperbolic Dehn filling of a hyperbolic $4$-manifold \cite{MR}.

%%%%%%%%%%%%%%%%%%%%%%%%%%%%%%%%%%%%%%%%%%%%%%%%%%%%%%%%%%%%%%%%%%%%%%%%%%%%%%%%%%%%%%%%%%%%%%%%%
\section{A right-angled $5$-polytope of finite volume}
In this section, we construct a hyperbolic right-angled $5$-polytope of finite volume concerned with our main theorem 
and study its geometry and combinatorics. 

\subsection{Construction of the right-angled $5$-polytope of finite volume}
We list the 48 unit space-like vectors in $\R^{5,1}$: 
\begin{itemize}
\item[(i)] 8 vectors obtained by permuting the first 4 coordinates of $\(\pm{1},0,0,0,1,1\)$. 
\item[(ii)] 8 vectors obtained by permuting the first 4 coordinates of $\(\pm{1},0,0,0,-\frac{1}{2},\frac{1}{2}\)$. 
\item[(iii)] 32 vectors obtained by permuting the first 4 coordinates of $\(\pm{1},\pm{1},\pm{1},0,\frac{1}{2}, \frac{3}{2}\)$.
\end{itemize}
For a vector $\bvec{v}$ of the items (i), (ii), and (iii), we denote the first 4 coordinates of $\bvec{v}$ by $\bvec{a}(\bvec{v})$. 
\begin{lem}\label{lem:3-1}
For two distinct vectors $\bvec{v}$ and $\bvec{w}$ of the items (i), (ii), and (iii), we have 
\begin{equation}
\<\bvec{v}, \bvec{w}\>=0, -1, -2, -3, -4, -5. \label{eq:3-(1)}
\end{equation}
Moreover, the $5$-polytope $\mathcal{P}$ defined by the 48 vectors of the items (i), (ii), and (iii) is right-angled. 
\end{lem}
\proof By the definition of the Lorentzian inner product $\<\bvec{v}, \bvec{w}\>$, we have 
\begin{equation}
\<\bvec{v}, \bvec{w}\>=\begin{cases}
         \(\bvec{a}(\bvec{v}), \bvec{a}(\bvec{w})\) &\text{ if }\bvec{v} \text{ and }\bvec{w} \text{ are of the item }(i). \\
         \(\bvec{a}(\bvec{v}), \bvec{a}(\bvec{w})\)-1 &\text{ if }\bvec{v}\text{ and }\bvec{w} \text{ is of the item }(i) \text{ and }(ii), \text{ respectively}. \\
         \(\bvec{a}(\bvec{v}), \bvec{a}(\bvec{w})\)-1 &\text{ if }\bvec{v}\text{ and }\bvec{w} \text{ is of the item }(i) \text{ and }(iii), \text{ respectively}. \\
         \(\bvec{a}(\bvec{v}), \bvec{a}(\bvec{w})\) &\text{ if }\bvec{v}\text{ and }\bvec{w} \text{ are of the item }(ii). \\
         \(\bvec{a}(\bvec{v}), \bvec{a}(\bvec{w})\)-1 &\text{ if }\bvec{v}\text{ and }\bvec{w} \text{ is of the item }(ii) \text{ and }(iii), \text{ respectively}. \\
         \(\bvec{a}(\bvec{v}), \bvec{a}(\bvec{w})\)-2 &\text{ if }\bvec{v}\text{ and }\bvec{w} \text{ are of the item }(iii).
         \end{cases}\label{eq:3-(2)}
\end{equation}
The equations \eqref{eq:3-(2)}, together with the fact that $\(\bvec{a}(\bvec{v}), \bvec{a}(\bvec{w})\)=1, 0, -1, -2, -3$, imply the equation \eqref{eq:3-(1)}. 

The fact that 2 facets of $\mathcal{P}$ defined by $\bvec{v}$ and $\bvec{w}$ intersects if and only if $\abs{\<\bvec{v}, \bvec{w}\>}<1$, 
together with the definition of the dihedral angles, implies that the $5$-polytope $\mathcal{P}$ is right-angled. \qed

\vspace{2mm}
Using the computer program CoxIter \cite{Gu} we have that the polytope $\mathcal{P}$ has finite volume . 
In this paper, 
our objective is to study local rigidity of the RACGs with Fuchsian ends 
obtained from the RACG $\Gamma_{\mathcal{P}}$ associated with $\mathcal{P}$. 

\vspace{2mm}
\subsection{Combinatorics of $\mathcal{P}$.}
\begin{figure}[htbp]
\centering
\includegraphics[scale=.75]{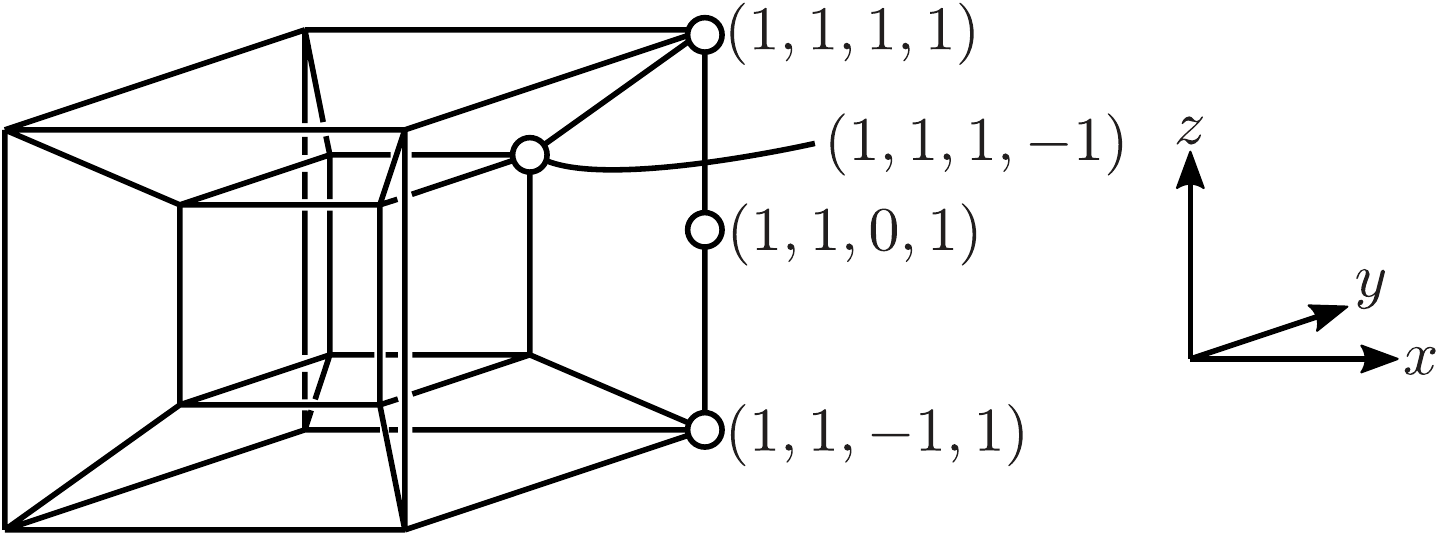}
\caption{The diagram $\mathcal{C}$. }
\label{fig:figure1}
\end{figure}

In order to see the mutual positions of two facets of $\mathcal{P}$, 
we use the diagram $\mathcal{C}$ as in Figure \ref{fig:figure1}. 
The $8$ cubes in Figure \ref{fig:figure1} correspond to the vectors 
$(\pm{1},0,0,0)$, $(0, \pm{1}, 0, 0)$, $(0,0,\pm{1},0)$, and $(0,0,0,\pm{1})$ in $\R^4$. 
For example, the right, back, upper, and outer cube corresponds to $(1,0,0,0), (0,1,0,0), (0,0,1,0)$, and $(0,0,0,1)$, respectively. 
Since every quadrilateral in Figure \ref{fig:figure1} is an intersection of exactly 2 cubes, 
we correspond the quadrilaterals to the sums of the 2 vectors corresponding to the cubes. 
Similarly, 
the edges and vertices in Figure \ref{fig:figure1} correspond to the sums of the 3 and 4 vectors corresponding to the cubes, 
respectively. 
We associate the facets of $\mathcal{P}$ with the cubes or edges of $\mathcal{C}$ according to $\bvec{a}(\bvec{v})$, 
where $\bvec{v}$ is the defining unit space-like vector of $\mathcal{P}$. 
Thus we get a map $\phi$ from the facets of $\mathcal{P}$ to the cubes and edges of $\mathcal{C}$:
for a facet $F_{\bvec{v}}$ of $\mathcal{P}$ defined by $\bvec{v}$, 
\begin{equation}
\phi(F_{\bvec{v}})=\text{ the cube or edge corresponding to }\bvec{a}(\bvec{v})\in{\R^4}. 
\end{equation}
Using the map $\phi$ we can read off the mutual position of two facets of $\mathcal{P}$ from $\mathcal{C}$. 
We list such a correspondence between the mutual positions of the facets of $\mathcal{P}$ and the combinatorics of $\mathcal{C}$. 
\begin{lem}\label{lem:3-2}
The followings hold for facets $F_{\bvec{v}}$ and $F_{\bvec{w}}$ of $\mathcal{P}$. 
\begin{enumerate}
\item[(i)] If $F_{\bvec{v}}$ and $F_{\bvec{w}}$ are of type (i), then they are 
\begin{equation*}
\begin{cases}
\text{ intersecting } &\Leftrightarrow \text{ the cubes } \phi(F_{\bvec{v}}) \text{ and }\phi(F_{\bvec{w}}) \text{ are intersecting}. \\
\text{ parallel }&\Leftrightarrow \text{ the cubes } \phi(F_{\bvec{v}}) \text{ and }\phi(F_{\bvec{w}}) \text{ are disjoint}. 
\end{cases}
\end{equation*}
\item[(ii)] If $F_{\bvec{v}}$ and $F_{\bvec{w}}$ is of type (i) and (ii), respectively, then they are 
\begin{equation*}
\begin{cases}
\text{ intersecting } &\Leftrightarrow \text{ the cubes }\phi(F_{\bvec{v}}) \text{ and }\phi(F_{\bvec{w}}) \text{ coincide with each other}. \\
\text{ parallel } & \Leftrightarrow \text{ the cubes }\phi(F_{\bvec{v}}) \text{ and }\phi(F_{\bvec{w}}) \text{ are intersecting }. \\
\text{ ultraparallel }&\Leftrightarrow \text{ the cubes } \phi(F_{\bvec{v}}) \text{ and }\phi(F_{\bvec{w}}) \text{ are disjoint}. 
\end{cases}
\end{equation*}
\item[(iii)] If $F_{\bvec{v}}$ and $F_{\bvec{w}}$ is of type (i) and (iii), respectively, then they are 
\begin{equation*}
\begin{cases}
\text{ intersecting } &\Leftrightarrow \text{ the edge }\phi(F_{\bvec{w}}) \text{ is contained in the cube }\phi(F_{\bvec{v}}) \\
\text{ parallel }&\Leftrightarrow \text{ the edge }\phi(F_{\bvec{w}}) \text{ meets the cube }\phi(F_{\bvec{v}})\text{ at a vertex}. \\
\text{ ultraparallel }&\Leftrightarrow \text{ the edge }\phi(F_{\bvec{w}}) \text{ is disjoint from the cube }\phi(F_{\bvec{v}}). 
\end{cases}
\end{equation*}
\item[(iv)] If $F_{\bvec{v}}$ and $F_{\bvec{w}}$ are of type (ii), then they are 
\begin{equation*}
\begin{cases}
\text{ intersecting } &\Leftrightarrow \text{ the cubes } \phi(F_{\bvec{v}}) \text{ and }\phi(F_{\bvec{w}}) \text{ are intersecting}. \\
\text{ parallel }&\Leftrightarrow \text{ the cubes } \phi(F_{\bvec{v}}) \text{ and }\phi(F_{\bvec{w}}) \text{ are disjoint}. 
\end{cases}
\end{equation*}
\item[(v)] If $F_{\bvec{v}}$ and $F_{\bvec{w}}$ is of type (ii) and (iii), respectively, then they are 
\begin{equation*}
\begin{cases}
\text{ intersecting } &\Leftrightarrow \text{ the edge }\phi(F_{\bvec{w}}) \text{ is contained in the cube }\phi(F_{\bvec{v}}) \\
\text{ parallel }&\Leftrightarrow \text{ the edge }\phi(F_{\bvec{w}}) \text{ meets the cube }\phi(F_{\bvec{v}})\text{ at a vertex}. \\
\text{ ultraparallel }&\Leftrightarrow \text{ the edge }\phi(F_{\bvec{w}}) \text{ is disjoint from the cube }\phi(F_{\bvec{v}}). 
\end{cases}
\end{equation*}
\item[(vi)] If $F_{\bvec{v}}$ and $F_{\bvec{w}}$ are of type (iii), then they are 
\begin{equation*}
\begin{cases}
\text{ intersecting } &\Leftrightarrow \text{ the edges }\phi(F_{\bvec{v}}) \text{ and }\phi(F_{\bvec{w}})\text{ meet at a vertex}. \\
\text{ parallel }&\Leftrightarrow \phi(F_{\bvec{v}}) \text{ and }\phi(F_{\bvec{w}})\text{ are opposite edges in a quadrilateral}. \\
\text{ ultraparallel }&\Leftrightarrow \text{ otherwise }.
\end{cases}
\end{equation*}
\end{enumerate}
\end{lem}

\vspace{2mm}
\subsection{The symmetry group of $\mathcal{P}$.}
A \textit{symmetry of $\mathcal{P}$} is an isometry of $\H^5$ which preserves the $5$-polytope $\mathcal{P}$. 
The set of all the symmetries of $\mathcal{P}$ is denoted by $\Sym{(\mathcal{P})}$. 
Note that $\Sym{(\mathcal{P})}$ acts on the set $\mathcal{F}(\mathcal{P})$ of the facets of $\mathcal{P}$ 
such that the mutual positions of the facets are preserved. 
Define the action of the symmetric group $\mathfrak{S}_4$ on $\H^5$ 
\begin{equation}
\sigma\(x_1, x_2, x_3, x_4, x_5, x_6\)=\(x_{\sigma^{-1}(1)}, x_{\sigma^{-1}(2)}, x_{\sigma^{-1}(3)}, x_{\sigma^{-1}(4)}, x_5, x_6\). 
\end{equation}
It is easy to see that $\mathfrak{S}_4$ is a subgroup of $\Sym{(\mathcal{P})}$. 
Let $\mathcal{F}_{(i)}, \mathcal{F}_{(ii)}$, and $\mathcal{F}_{(iii)}$ denote the set of the facets of type (i), (ii), and (iii), respectively. 
Then, $\mathcal{S}_4$ acts on $\mathcal{F}_{(i)}$ and $\mathcal{F}_{(ii)}$ transitively. 
\begin{lem}\label{lem:3-3}
The followings hold for every symmetry $f$ of $\mathcal{P}$. 

(i) $f$ preserves $\mathcal{F}_{(i)}\sqcup{\mathcal{F}_{(ii)}}$ and $\mathcal{F}_{(iii)}$. 

(ii) If $f\(\mathcal{F}_{(i)}\)\cap{\mathcal{F}_{(ii)}}\neq{\emptyset}$, then $f\(\mathcal{F}_{(i)}\)=\mathcal{F}_{(ii)}$. 
\end{lem}
\proof (i) By Lemma \ref{lem:3-2}, we see that 
\begin{equation}
\text{ the number of }3\text{-faces of a facet }F\text{ of }\mathcal{P}=
\begin{cases}
24 & F\in{\mathcal{F}_{(i)}\sqcup{\mathcal{F}_{(ii)}}}. \\
10 & F\in{\mathcal{F}_{(iii)}}. 
\end{cases}
\end{equation}
Therefore every symmetry $f$ of $\mathcal{P}$ preserves $\mathcal{F}_{(i)}\sqcup{\mathcal{F}_{(ii)}}$ and $\mathcal{F}_{(iii)}$.

\vspace{1mm}
(ii) Suppose that $f$ maps a facet $F_1\in{\mathcal{F}_{(i)}}$ to a facet $F_2\in{\mathcal{F}_{(ii)}}$. 
We can assume, 
by the action of $\mathfrak{S}_4$ on $\mathcal{F}_{(i)}$ and $\mathcal{F}_{(ii)}$ if necessary, 
that the defining vectors of $F_1$ and $F_2$ are $(1,0,0,0,1,1)$ and $(1,0,0,0,-1/2,1/2)$. 
If a facet $F$ intersects $F_1$, then the facet $f(F)$ intersects $F_2$. 
By Lemma \ref{lem:3-2}, 
we see that $f$ maps every facet $F$ of type (i) intersecting $F_1$ to the facet $f(F)$ of type (ii) intersecting $F_2$. 
Thus we have that $f\(\mathcal{F}_{(i)}\)=\mathcal{F}_{(ii)}$. \qed

\vspace{2mm}
A bijection on the set of the cubes of $\mathcal{C}$ preserving its mutual positions is called a \textit{symmetry} of $\mathcal{C}$. 
We denote by $\Sym{(\mathcal{C})}$ the set of all the symmetries of $\mathcal{C}$. 
It is known that $\Sym{(\mathcal{C})}$ is isomorphic to the finite Coxeter group $B_4$. 

\begin{lem}
Every symmetry $f$ of $\mathcal{P}$ induces the symmetry $\phi_\ast f$ of $\mathcal{C}$. 
Moreover, the map $\phi_\ast:\Sym{(\mathcal{P})}\to{\Sym{(\mathcal{C})}}$ is a group homomorphism. 
\end{lem}
\proof Set $A=\{(\pm{1},0,0,0), (0,\pm{1},0,0),(0,0,\pm{1},0),(0,0,0,\pm{1})\}\subset{\R^4}$. 
For every element $\bvec{e}\in{A}$, define $\underline{f}(\bvec{e})=\bvec{a}\(f(\bvec{v})\)$, where $\bvec{v}$ is a defining vector of $\mathcal{P}$ satisfying $\bvec{a}(\bvec{v})=\bvec{e}$. 
By Lemma \ref{lem:3-2} and \ref{lem:3-3}, we see that $\underline{f}(\bvec{e})$ does not depend on the choice of $\bvec{v}$. 
The following diagram commutes. 
\[
\xymatrix{
\mathcal{F}_{(i)}\sqcup{\mathcal{F}_{(ii)}} \ar[r]^-{a} \ar[d]_-{f} & A \ar[d]^-{\underline{f}} \\
\mathcal{F}_{(i)}\sqcup{\mathcal{F}_{(ii)}} \ar[r]_-{a} & A
}
\]
The set $A$ is canonically identified with the set of the cubes of $\mathcal{C}$. 
Therefore we get the symmetry $\phi_\ast f$ of $\mathcal{C}$ induced by $f$. 
It is easy to see that the map $\phi_\ast:\Sym{(\mathcal{P})}\to{\Sym{(\mathcal{C})}}$ is a group homomorphism. \qed

\begin{cor}
The symmetry group $\Sym{(\mathcal{P})}$ is generated by 5 reflections $r_1, \cdots, r_5$ associated with the following vectors: 
\begin{gather*}
 \bvec{u}_1=(1,0,0,0,0,0), \quad
 \bvec{u}_2=\(-\frac{1}{\sqrt{2}},\frac{1}{\sqrt{2}},0,0,0,0\), \quad
 \bvec{u}_3=\(0,-\frac{1}{\sqrt{2}}, \frac{1}{\sqrt{2}},0,0,0\), \\
 \bvec{u}_4=\(0,0,-\frac{1}{\sqrt{2}},\frac{1}{\sqrt{2}},0,0\), \quad
 \bvec{u}_5=\(0,0,0,0,-\frac{3}{2\sqrt{2}},-\frac{1}{2\sqrt{2}}\).
\end{gather*}
\end{cor}
\proof Since every $f\in{\ker{\phi_\ast}}$ acts on $A$ as $id_A$, we have that $\ker{\phi_\ast}=\{id_{\H^5}, r_5\}$. 
The fact that the reflections $r_i \ (i=1, \cdots, 4)$ and $r_5$ are commutative 
implies that the group homomorphism $\phi_\ast$ is injective on the subgroup $\<r_1, \cdots, r_4\>$. 
The group homomorphism $\phi_\ast$ is surjective, 
since the subgroup $\<r_1, \cdots, r_4\>$ and $\Sym{(\mathcal{C})}$ is isomorphic to the finite Coxeter group $B_4$. 
Thus $\Sym{(\mathcal{P})}$ is generated by the reflections $r_1, \cdots, r_5$. \qed

\section{Local rigidity of $\mathcal{P}$}
We label the 48 vectors of the items (i), (ii), and (iii) as follows. 
\begin{align*}
\text{(i) : }&X_{\pm}=(\pm{1},0,0,0,1,1), \quad Y_{\pm}=(0,\pm{1},0,0,1,1), \\
           &Z_{\pm}=(0,0,\pm{1},0,1,1), \quad W_{\pm}=(0,0,0,\pm{1},1,1). \\
\text{(ii) : }&S_{X_{\pm}}=\(\pm{1},0,0,0,-\frac{1}{2}, \frac{1}{2}\), \quad S_{Y_{\pm}}=\(0,\pm{1},0,0,-\frac{1}{2}, \frac{1}{2}\), \\
            &S_{Z_{\pm}}=\(0,0,\pm{1},0,-\frac{1}{2}, \frac{1}{2}\), \quad S_{W_{\pm}}=\(0,0,0,\pm{1},-\frac{1}{2}, \frac{1}{2}\). \\
\text{(iii) : }& S_{\bvec{a}}=\(\bvec{a}, \frac{1}{2}, \frac{3}{2}\)\text{ where }\bvec{a}\in{\R^4}\text{ is obtained from permuting}\\
&\text{ the coordinates of }(\pm{1}, \pm{1}, \pm{1}, 0). 
\end{align*}
For abbreviation, we use the same letter $A$ (resp. $A^-$) for the hyperplane $H_A$ (resp. the closed half-space $H_A^-$) of $\H^d$ with the outer normal vector $A$, where $A$ is a vector labeled above. 

\vspace{2mm}
Suppose that two sets $\mathcal{F}$ and $\mathcal{F}'$ of disjoint facets of $\mathcal{P}$ satisfy that 
$\mathcal{F}'\subset{\mathcal{F}}$. 
Then, if the group $\Gamma_{\mathcal{F}}$ is locally rigid, so is the group $\Gamma_{\mathcal{F}'}$. 
The observation implies that in order to prove that the RACGs obtained from $\Gamma_{\mathcal{P}}$ are locally rigid, 
it is sufficient to consider the case that 
$\mathcal{F}$ is an inclusion-wise maximal element among all the sets of disjoint facets of $\mathcal{P}$. 
From now on, $\mathcal{F}$ denotes such a set of disjoint facets of $\mathcal{P}$. 

\begin{lem}\label{lem:local rigidity}
For every $\rho\in{U_0}$, there exists an isometry $f\in{\Isom{(\H^5)}}$ such that $f(A_{\pm}(\rho))=A_{\pm}(\rho_0)$ and $f(S_{A_{\pm}}(\rho))=S_{A_{\pm}}(\rho_0)$ where $A=X, Y, Z, W$. 
\end{lem}
\proof Since $\Sym{(\mathcal{P})}$ acts on the generating set $S$ by permutation, 
it is sufficient to show the assertion for the case that $\mathcal{F}$ contains $W_{\pm}$ and $S_{W_{\pm}}$. 
Throughout the proof, we use the upper half-space model $\U^5$ of hyperbolic $5$-space. 
Consider the subgroup $\Gamma_{\text{cube}}$ of $\Gamma_{\mathcal{F}}$ generated by the 6 reflections $\rho(r_{X_{\pm}}), \rho(r_{Y_{\pm}})$, and $\rho(r_{Z_{\pm}})$. 
Since $\Gamma_{\text{cube}}$ fixes a unique point at infinity, 
by an isometry $f_0\in{\Isom{(\U^5)}}$, we assume that the point at the infinity is $\infty\in{\partial{\U^5}}$. 
Then the hyperplanes $A_{\pm}(\rho)$ are Euclidean hyperplanes and the hyperplanes $S_{A_{\pm}}(\rho)$ are Euclidean hemispheres orthogonal to $\R^4\subset{\partial{\U^5}}$. 
By abuse of notation, 
we continue to write $A_{\pm}$ and $S_{A_{\pm}}$ for the boundaries of the hyperplanes in $\overline{\U}^5$. 
Since the stabilizer subgroup of $\Isom{(\U^5)}$ with respect to $\infty$ is the similarity group $\Sim{(\R^4)}$ by the Poincar\'e extension, 
there exists $f_1\in{\Sim{(\R^4)}}$ such that 
\begin{align}
f_1(X_-^-(\rho))&=\Set{(x,y,z,w)\in{\R^4}|x\geq{0}}. 
\end{align}
The outer normal vector of $f_1(X_-^-(\rho))$ is $(-1,0,0,0)$. 
Thus we can see that every similarity $f\in{\Sim{(\R^4)}}$ which preserves the closed half-space $f_1(X_-^-)$ is written as the following form. 
\begin{equation}
f(\bvec{x})=\lambda\begin{pmatrix} 1 & \bvec{0} \\ \bvec{0} & A\end{pmatrix}\bvec{x}+\begin{pmatrix} 0 \\ \bvec{b} \end{pmatrix}
\text{ where }\lambda\in{\R_{>0}}, \ A\in{O(3)}\text{, and }\bvec{b}\in{\R^3}. \label{eq:(7)}
\end{equation}
The relation that $\(\rho(r_{X_-})\rho(r_{Y_-})\)^2=id_{\U^5}$ implies that the hyperplanes $X_-(\rho)$ and $Y_-(\rho)$ intersect orthogonally. 
Therefore the closed half-space $f_1(Y_-^-)$ can be written as 
\begin{equation}
f_1(Y_-^-(\rho))=\Set{\bvec{x}\in{\R^4}|(\bvec{x}, (0,y,z,w))\leq{t}}. 
\end{equation}
Since the outer normal vector of $f_1(Y_-^-)$ is of the form $(0, y, z, w)$, there exists an isometry $f_2\in{\Isom{(\R^4)}}$ of the form \eqref{eq:(7)} such that 
\begin{align}
f_2f_1(X_-^-(\rho))&=\Set{(x, y, z, w)\in{\R^4}|x\geq{0}}. \\
f_2f_1(Y_-^-(\rho))&=\Set{(x, y, z, w)\in{\R^4}|y\geq{0}}. 
\end{align}
Since $f_2f_1(S_{X_-}(\rho))$ intersects with $f_2f_1(X_-(\rho))$ orthogonally and parallel to $f_2f_1(Y_-(\rho))$, 
we obtain that 
\begin{equation}
f_2f_1(S_{X_-}(\rho))=\Set{(x,y,z,w)\in{\R^4}|\abs{(0,a,b,c)-(x,y,z,w)}^2=a^2} \ (a>0, \ b,c\in{\R}). 
\end{equation}
By $f_3(\bvec{x})=\frac{1}{a}(\bvec{x}-(0,0,b,c))$, 
\begin{align}
f_3f_2f_1(X_-^-(\rho))&=\Set{(x, y, z, w)\in{\R^4}|x\geq{0}}. \\
f_3f_2f_1(Y_-^-(\rho))&=\Set{(x, y, z, w)\in{\R^4}|y\geq{0}}. \\
f_3f_2f_1(S_{X_-}(\rho))&=\Set{(x,y,z,w)\in{\R^4}|\abs{(0,1,0,0)-(x,y,z,w)}^2=1}. 
\end{align}
Similarly, by considering the facets intersecting with $f_3f_2f_1(S_{Y_-}(\rho))$ orthogonally, we have that 
\begin{equation}
f_3f_2f_1(S_{X_-}(\rho))=\Set{(x,y,z,w)\in{\R^4}|\abs{(d,0,0,0)-(x,y,z,w)}^2=d^2} \ (d>0). 
\end{equation}
We can write every symmetry $f\in{\Sim{(\R^4)}}$ which preserves $X_-(\rho), Y_-(\rho), S_{X_-}(\rho)$, and $S_{Y_-}(\rho)$ in the form 
\begin{equation}
f(\bvec{x})=\begin{pmatrix} I_2 & \bvec{0} \\ \bvec{0} & B\end{pmatrix}\bvec{x} \text{ where }B\in{O(2)}. \label{eq:(16)}
\end{equation}
Since $f_3f_2f_1(Z_-(\rho))$ intersects with $f_3f_2f_1(X_-(\rho))$ and $f_3f_2f_1(Y_-(\rho))$ and parallel to $f_3f_2f_1(S_{X_-}(\rho))$ and $f_3f_2f_1(S_{Y_-}(\rho))$, 
by an isometry $f_4$ of the form \eqref{eq:(16)}, 
\begin{equation}
f_4f_3f_2f_1(Z_-(\rho))=\Set{(x,y,z,w)\in{\R^4}|z\geq{-1}} \text{ and }d=1. 
\end{equation}
Let us denote by $f$ the isometry $f_5f_4f_3f_2f_1\in{\Isom{(\U^5)}}$ where $f_5(\bvec{x})=\bvec{x}-(1,1,0,0)$. 
Then, 
\begin{align}
f(X_-^-(\rho))&=\Set{(x,y,z,w)|x\geq{-1}}=X_-^-(\rho_0). \\
f(Y_-^-(\rho))&=\Set{(x,y,z,w)|y\geq{-1}}=Y_-^-(\rho_0). \\
f(Z_-^-(\rho))&=\Set{(x,y,z,w)|z\geq{-1}}=Z_-^-(\rho_0). \\
f(S_{X_-}(\rho))&=\Set{(x,y,z,w)|\abs{(x,y,z,w)-(-1,0,0,0)}^2=1}=S_{X_-}(\rho_0). \\
f(S_{Y_-}(\rho))&=\Set{(x,y,z,w)|\abs{(x,y,z,w)-(0,-1,0,0)}^2=1}=S_{Y_-}(\rho_0). 
\end{align}
The facts that the hyperplane $X_+(\rho)$ is parallel to $X_-(\rho)$ and intersects with $S_{Y_-}(\rho)$ imply that 
\begin{equation}
f(X_+^-(\rho))=\Set{(x,y,z,w)\in{\R^4}|x\leq{1}}=X_+^-(\rho_0). 
\end{equation}
Similarly, we have that $f(A_{\pm}(\rho))=A_{\pm}(\rho_0)$ and $f(S_{A_{\pm}}(\rho))=S_{A_{\pm}}(\rho_0)$. \qed

\begin{theo}
For every $\rho\in{U_0}$, 
there exists an isometry $f\in{\Isom{(\H^5)}}$ such that $f\rho f^{-1}=\rho_0$, that is, 
all the RACGs with Fuchsian ends obtained from $\Gamma_{\mathcal{P}}$ are locally rigid. 
\end{theo}
\proof By Lemma \ref{lem:local rigidity}, 
there exists $f\in{\Isom{(\U^5)}}$ such that $f(A_{\pm}(\rho))=A_{\pm}(\rho_0)$ and $f(S_{A_{\pm}}(\rho))=S_{A_{\pm}}(\rho_0)$. 
For simplicity of notation, we write $r_H^\rho$ and $r_H^0$ for $\rho(r_H)$ and $\rho_0(r_H)$, respectively. 
Since $fr_Hf^{-1}=r_{f(H)}$, we have that 
\begin{equation}
fr^\rho_{A_{\pm}}f^{-1}=r_{A_{\pm}}^0 \text{ and }fr^\rho_{S_{A_{\pm}}}f^{-1}=r_{S_{A_{\pm}}}^0, 
\end{equation}
where $A\in{\{X, Y, Z, W\}}$. 
Consider a facet $S_{\bvec{a}}\not\in{\mathcal{F}}$. 
Then we show that $f(S_{\bvec{a}}(\rho))=S_{\bvec{a}}(\rho_0)$. 
For an example, we shall show the equality for the case that $\bvec{a}=(1,1,1,0)$. 
Since the edge of $\mathcal{C}$ corresponding to $\bvec{a}=(1,1,1,0)$ is the intersection of 3 cubes 
corresponding to $(1,0,0,0), (0,1,0,0)$, and $(0,0,1,0)$ (see Figure \ref{fig:figure1}), 
the number of elements of $\mathcal{F}$ of type (i) and (ii) corresponding to one of the 3 cubes is at most one, respectively. 
%at most one hyperplane of type (i) and (ii) corresponding to one of the 3 cubes is an element of $\mathcal{F}$, respectively. 
We divide the proof into 2 cases; 
one is for the case that such a hyperplane does not exist, and the other is converse. 
If such a hyperplane does not exist, 
the fact that the hyperplanes $X_{+}, Y_+, S_{X_+}, S_{Y_+}, Z_+$ intersect with $S_{\bvec{a}}$ orthogonally 
implies that 
\begin{equation}
S_{\bvec{a}}=\Set{(x,y,z,w)\in{\R^4}|\abs{(x,y,z,w)-(1,1,1,0)}^2=1}. 
\end{equation}
If exists, we may assume that such a hyperplane is $Z_+$ and $S_{Z_+}$ according to its types. 
Since $S_{\bvec{a}}$ intersect with the hyperplanes $X_{+}, Y_+, S_{X_+}, S_{Y_+}$ orthogonally and parallel to $W_+$, 
we have that 
\begin{equation}
S_{\bvec{a}}=\Set{(x,y,z,w)\in{\R^4}|\abs{(x,y,z,w)-(1,1,1,0)}^2=1}. 
\end{equation}
Therefore we see that $f(S_{\bvec{a}}(\rho))=S_{\bvec{a}}(\rho_0)$. \qed

\section{Acknowledgements}
The author wishes to express his gratitude to Stefano Riolo for having a time for fruitful discussions.

\end{document}